\documentclass[11pt]{article}
\usepackage{amssymb, amscd}
\usepackage{amsfonts}
\usepackage{amsmath}
\usepackage{graphicx}

\newtheorem{thm}{Theorem}

\newtheorem{cor}[thm]{Corollary}

\newtheorem{defn}[thm]{Definition}
\newtheorem{exmp}[thm]{Example}

\newtheorem{prop}[thm]{Proposition}
\newtheorem{rem}[thm]{Remark}

\begin{document}

\title{\bf Parallel transport of $Hom$-complexes\\ and the Lov\' asz conjecture}
\author{Rade T. \v Zivaljevi\' c }

\date{May 15, 2005}

\maketitle

\begin{abstract}
 The groupoid of projectivities, introduced by M.~Joswig
 \cite{Josw2001}, serves as a basis for a construction of parallel
 transport of graph and more general $Hom$-complexes. In this
 framework we develop a general conceptual approach to the Lov\'{a}sz
 conjecture, recently resolved by E.~Babson and D.~Kozlov in
 \cite{BabsonKozlov2}, and extend their result from graphs to the
 case of simplicial complexes.
\end{abstract}

\section{Introduction\label{sec_introduction}}

Hidden in the background of the Babson and Kozlov proof of Lov\'
asz conjecture \cite{BabsonKozlov2} \cite{Kozlov-review} are
interesting topological and combinatorial concepts and structures
associated to graphs and graph homomorphisms. The proof itself
runs in two phases, each phase divided in several steps, often
involving a detailed case analysis. For this reason the underlying
secondary structures may not be visible or immediately recognized
under the layers of intricate technical details. Recall that the
crux of Babson and Kozlov approach is a skilful and technically
quite involved application of spectral sequences. One of the
classical applications of this method is to the (co)homology of
fibered spaces which by the nature are spaces which alow some form
of transport from one fibre to another.

In this paper we focus on one of these secondary structures which
can be, somewhat informally, described as the ``parallel
transport'' of graph complexes over graphs.

The introduction of this structure and recognition of its role
leads to a great simplification of the proof of Lov\' asz
conjecture in some cases. Another of its features, aside from
offering a conceptual ``explanation'' for the success of one of
the technical approaches of Babson and Kozlov, is its potential
for generating other statements of this type.

The ``parallel transport'' of graph complexes introduced here
seems to be a novel concept. However the groupoids (groups of
projectivities) used in its definition have already appeared in
geometric combinatorics in the work of M.~Joswig \cite{Josw2001},
see also \cite{Izm2001} \cite{Izm-Usp2001} \cite{IzmJos2002},
where they have been applied to toric manifolds, branched
coverings over $S^3$, colorings of simple polytopes etc.

It is an exciting ``coincidence'' that there have been other
recent developments in geometric and algebraic combinatorics where
groupoids and associated objects and constructions were {\em
implicitly} used, \cite{BBLL} \cite{BL} \cite{BKLW}
\cite{Gaif2004}. These are not isolated examples of course. In
particular one should be fully aware of a rich and deep
combinatorics already present in numerous categorical
constructions related to groupoids and their applications in
geometry and mathematical physics.

\section{The Lov\' asz conjecture}

One of central themes in topological combinatorics, after the
landmark paper of Laszlo Lov\' asz \cite{Lovasz} where he proved
the classical Kneser conjecture, has been the study and
applications of graph complexes.

The underlying theme is to explore how the topological complexity
of a graph complex $X(G)$ reflects in the combinatorial complexity
of the graph $G$ itself. The results one is usually interested in
come in the form of inequalities $\alpha(X(G))\leq \xi(G)$, or
equivalently in the form of implications
\[
\alpha(X(G))\geq p \Rightarrow \xi(G)\geq q,
\]
where $\alpha(X(G))$ is a topological invariant of $X(G)$, while
$\xi(G)$ is a combinatorial invariant of the graph $G$.

The most interesting candidate for the invariant $\xi$ has been
the chromatic number $\chi(G)$ of $G$, while the role of the
invariant $\alpha$ was played by the ``connectedness'' of $X(G)$,
its equivariant index, the height of an associated characteristic
cohomology class etc., see \cite{Kozlov-review} \cite{Matousek}
\cite{MatZieg} \cite{Ziv2005} for recent accounts.

The famous result of Lov\' asz quoted above is today usually
formulated in the form of an implication
\begin{equation}
Hom(K_2,G) \mbox{ {\rm is $k$-connected} } \Rightarrow \chi(G)\geq
k+3,
 \end{equation}
where $Hom(K_2,G)$ is the so called ``box complex'' of $G$. The
box complex is a special case of a general graph complex
$Hom(H,G)$ (also introduced by L.~Lov\' asz), a cell complex which
functorially depends on the input graphs $H$ and $G$.

An outstanding conjecture in this area, refereed to as ``Lov\' asz
conjecture'', was that one obtains a better bound if the graph
$K_2$ in (1) is replaced by an odd cycle $C_{2r+1}$. More
precisely Lov\' asz conjectured that
\begin{equation}\label{eq:Lovasz-conj}
Hom(C_{2r+1},G) \mbox{ {\rm is $k$-connected} } \Rightarrow
\chi(G)\geq k+4.
 \end{equation}
This conjecture was confirmed by Babson and Kozlov in
\cite{BabsonKozlov2}, see also \cite{Kozlov-review} for a more
detailed exposition.

Our objective is to develop methods which both offer a simplified
approach to the proof of implication (\ref{eq:Lovasz-conj}), at
least in the case when $k$ is odd, and providing new insight, open
a possibility of proving similar results for other classes of
(hyper)graphs and simplicial complexes.

An example of such a result is Theorem~\ref{thm:main}. One of its
corollaries is the following implication,
\begin{equation}\label{eq:complex-conn}
Hom(\Gamma,K) \mbox{ {\rm is $k$-connected} } \Rightarrow
\chi(K)\geq k+d+3
\end{equation}
which, under a suitable assumption on the complex $\Gamma$ and the
assumption that integer $k$ is odd, extends (\ref{eq:Lovasz-conj})
to the case of pure $d$-dimensional simplicial complexes.

\section{Parallel transport of $Hom$-complexes}

\subsection{Generalities about ``parallel transport''}
\label{sec:generalities}

In order to avoid any ambiguities, we briefly clarify what is in
this paper meant by a ``parallel transport'' on a ``bundle'' of
spaces.

A ``bundle'' is a map $\phi : X\rightarrow S$. We assume that $S$
is a set and that $X(i):=\phi^{-1}(i)$ is a topological space, so
a bundle is just a collection of spaces (fibres) $X(i)$
parameterized by $S$. If all spaces $X(i)$ are homeomorphic to a
fixed ``model'' space, this space is referred to as the fiber of
the bundle $\phi$.

Suppose that $\mathcal{G}$ is a groupoid on $S$ as the set of
objects. In other words $\mathcal{G}=(Ob(\mathcal{G}),
Mor(\mathcal{G}))$ is a small category where $Ob(\mathcal{G})=S$,
such that all morphisms $\alpha\in Mor(\mathcal{G})$ are
invertible.

A ``connection'' or ``parallel transport'' on the bundle
$\mathcal{X}=\{X(i)\}_{i\in S}$ is a functor $\mathcal{P}:
\mathcal{G}\rightarrow Top$ such $X(i)=\mathcal{P}(i)$ for each
$i\in S$.

Informally speaking, the groupoid $\mathcal{G}$ provides a ``road
map'' on $S$, while the functor $\mathcal{P}$ defines the
associated transport from one fibre to another.

Sometimes it is convenient to view the bundle
$\mathcal{X}=\{X(i)\}_{i\in S}$ as a map $\mathcal{X}:
S\rightarrow Top$. Then to define a ``connection'' on this bundle
is equivalent to enriching the map $\mathcal{X}$ to a functor
$\mathcal{P}: \mathcal{G}\rightarrow Top$.

\subsection{Natural bundles and groupoids over simplicial complexes}
\label{sec:natural}

Suppose that $K$ and $L$ are finite simplicial complexes and let
$k$ be an integer such that $0\leq k\leq {\rm dim}(K)$. Let $S_k =
S_k(K)$ be the set of all $k$-dimensional simplices in $K$. Define
a bundle $\mathcal{F}^L_k : S_k\rightarrow Top$ by the formula
\begin{equation}
\mathcal{F}^L_k(\sigma)=Hom(\sigma, L)\cong L_\Delta^{k+1}
\end{equation}
where $Hom(\sigma, L)$ is one of the $Hom$-complexes introduced in
Section~\ref{sec:Hom} and $L_\Delta^{k+1}$ is the complex well
known in topological combinatorics as the deleted product of $L$,
\cite{Matousek} Chapter~6. A typical cell in $L_\Delta^{k+1}$ is
of the form $e=\sigma_0\times\sigma_1\times\ldots\times\sigma_k\in
L^{k+1}$ where $\{\sigma_i\}_{i=0}^k$ is a collection of non-empty
simplices in $L$ such that if $i\neq j$ then
$\sigma_i\cap\sigma_j=\emptyset$. The corresponding cell in
$Hom(\sigma, L)$ is a function $\eta : V(\sigma)\rightarrow
L\setminus\{\emptyset\}$, where $V(\sigma)=\sigma^{(0)}$ is the
set of all vertices of $\sigma$, and if $v_1\neq v_2$ then
$\eta(v_1)\cap\eta(v_2)=\emptyset$.

\medskip\noindent
{\bf Example:} It is well known that if $L\cong \sigma^m =
\Delta^{[m+1]}$ is a $m$-dimensional simplex, then the associated
deleted square $(\sigma^m)^2_\Delta$ is homeomorphic to a
$(m-1)$-dimensional sphere. In other words,
$\mathcal{F}^{\sigma^m}_1 : S_1(K)\rightarrow Top$ is a spherical
bundle naturally associated to the simplicial complex $K$.

\medskip
Our next goal, in the spirit of Section~\ref{sec:generalities}, is
to identify a groupoid on the set $S_k$ which acts on the bundle
$\mathcal{F}^L_k$, i.e.\ a groupoid which provides a parallel
transport of fibres of the bundle $\mathcal{F}^L_k$. It is a
pleasant coincidence that this groupoid has already appeared in
geometric combinatorics \cite{Jos-Usp} \cite{Josw2001}. Indeed,
the {\em groups of projectivities} M.~Joswig introduced and
studied in these papers are just the {\em vertex} or {\em
isotropy} groups of a groupoid which we call the $k$-th {\em
groupoid of projectivities} of $K$ and denote by
$\mathcal{G}_k^P(K)$. In these and in subsequent papers
\cite{Izm2001} \cite{Izm-Usp2001} \cite{IzmJos2002}, the groups of
projectivities found interesting applications to toric manifolds,
branched coverings over $S^3$, colorings of simple polytopes, etc.
Here is a summary of this construction.

Two $k$-dimensional simplices $\sigma_0$ and $\sigma_1$ in $K$ are
called adjacent if they share a common $(k-1)$-dimensional face
$\tau$. A {\em perspectivity} from $\sigma_0$ to $\sigma_1$ is the
unique non-degenerated simplicial map
$\overrightarrow{\sigma_0\sigma_1} =
\langle\sigma_0,\sigma_1\rangle : \sigma_0\rightarrow\sigma_1$
which leaves the simplex $\tau$ point-wise fixed. In the special
case when $\sigma_0=\sigma_1$, the perspectivity
$\langle\sigma_0,\sigma_0\rangle : \sigma_0\rightarrow\sigma_0$ is
the identity map $I_{\sigma_0}$.

A {\em projectivity} between two, not necessarily adjacent,
simplices $\sigma_0$ and $\sigma_n$ is a composition of
perspectivities
\[
\langle\mathfrak{p}\rangle = \overrightarrow{\sigma_0\sigma_1}\ast
\overrightarrow{\sigma_1\sigma_2}\ast\ldots\ast
\overrightarrow{\sigma_{n-1}\sigma_n}
\]
where $\mathfrak{p}=(\sigma_0,\sigma_1,\ldots,\sigma_n)$ is a path
of $k$-dimensional simplices in $K$ such that $\sigma_{i-1}$ and
$\sigma_i$ share a common $(k-1)$-dimensional face $\tau_i$.

\medskip\noindent
{\bf Caveat:} Here we adopt a useful convention that $(x)(f\ast
g)= (g\circ f)(x)$ for each two composable maps $f$ and $g$. The
notation $f\ast g$ is often given priority over the usual $g\circ
f$ if we want to emphasize that the functions act on the points
from the right, that is if the arrows in the associated formulas
point from left to the right.

\begin{defn}{\rm \cite{Jos-Usp} \cite{Josw2001}}
The $k$-th {\em groupoid of projectivities} $\mathcal{G}_k^P(K)$
of a simplicial complex $K$, or the $P_k$-groupoid associated to
$K$, is the small category
 \[\mathcal{G}_k^P(K) =
(Ob(\mathcal{G}_k^P(K)), Mor(\mathcal{G}_k^P(K)))
 \]
which has the set $S_k = Ob(\mathcal{G}_k^P(K))$ of all
$k$-dimensional simplices for the set of objects, and for each two
simplices $\sigma_0,\sigma_1\in S_k$, the associated morphism set
$Mor_{\mathcal{G}_k^P(K)}(\sigma_0,\sigma_1)$ is the collection of
all projectivities from $\sigma_0$ to $\sigma_1$. The associated
point (isotropy) groups
\[
\Pi_k(K,\sigma_0):= Mor_{\mathcal{G}_k^P(K)}(\sigma_0,\sigma_0)
\]
are called the {\em groups of projectivities} or the {\em
combinatorial holonomy groups} of $K$ .
\end{defn}

\begin{prop}
For each finite simplicial complex $K$ and an auxiliary
``coefficient'' complex $L$, there exists a canonical connection
$\mathcal{P}^L = \mathcal{P}^L_{K,k}$ on the bundle
$\mathcal{F}^L_k$. In other words the function $\mathcal{F}^L_k:
S_k\rightarrow Top$ can be enriched (extended) to a functor
 \[\mathcal{F}^L_k: S_k\rightarrow Top.\]
\end{prop}

\medskip\noindent
{\bf Proof:} If $\overrightarrow{\sigma_0\sigma_1}$ is a
perspectivity from $\sigma_0$ to $\sigma_1$ and if $\eta :
V(\sigma_1) \rightarrow 2^{V(L)}\setminus\{\emptyset\}$ is a cell
in $Hom(\sigma,L)$, then $\mathcal{P}^L : \mathcal{F}^L(\sigma_1)
\rightarrow \mathcal{F}^L(\sigma_0)$ is the map defined by
$\mathcal{P}^L(\overrightarrow{\sigma_0\sigma_1})(\eta):=
\overrightarrow{\sigma_0\sigma_1}\ast\eta$. More generally, if
$\langle\mathfrak{p}\rangle =
\overrightarrow{\sigma_0\sigma_1}\ast
\overrightarrow{\sigma_1\sigma_2}\ast\ldots\ast
\overrightarrow{\sigma_{n-1}\sigma_n}$ is a projectivity between
$\sigma_0$ and $\sigma_n$, then
\begin{equation}
\mathcal{P}^L(\langle\mathfrak{p}\rangle) =
\mathcal{P}^L(\overrightarrow{\sigma_0\sigma_1})\ast
\mathcal{P}^L(\overrightarrow{\sigma_1\sigma_2})\ast\ldots\ast
\mathcal{P}^L(\overrightarrow{\sigma_{n-1}\sigma_n})
\end{equation}
or in other words
\begin{equation}
\mathcal{P}^L(\langle\mathfrak{p}\rangle)(\eta) =
\overrightarrow{\sigma_0\sigma_1}\ast
\overrightarrow{\sigma_1\sigma_2}\ast\ldots\ast
\overrightarrow{\sigma_{n-1}\sigma_n}\ast\eta.
\end{equation}
It is clear from the construction that the map
$\mathcal{P}^L(\langle\mathfrak{p}\rangle)$ depends only on the
projectivity $\langle\mathfrak{p}\rangle$ and not on the
associated path $\mathfrak{p}$. \hfill$\square$

\subsection{Parallel transport of graph complexes}
\label{sec:functor}

The main motivation for introducing the parallel transport of
$Hom$-complexes is the Lov\'{a}sz conjecture and its
ramifications. This is the reason why the case of graphs and the
graph complexes deserves a special attention. Additional
justification for emphasizing graphs comes from the fact that
graph complexes $Hom(G,H)$ have been studied in numerous papers
and today form a well established part of graph theory and
topological combinatorics. The situation with simplicial complexes
is quite the opposite. In order to extend the theory of
$Hom$-complexes from graphs to the category of simplicial
complexes, many concepts should be generalized and the
corresponding facts established in a more general setting. One is
supposed to recognize the main driving forces and to isolate the
most desirable features of the theory. A result should be a
dictionary/glossary of associated concepts, cf.\ Table~1.
Consequently, Section~\ref{sec:functor} should be viewed as an
important preliminary step, leading to the more general theory
developed in Sections~\ref{sec:ramifications} and \ref{sec:main}.

In order to simplify the exposition we assume, without a serious
loss of generality, that all graphs $G=(V(G),E(G))$ are without
loops and multiple edges. In short, graphs are $1$-dimensional
simplicial complexes. Let $G_{\overline{xy}}\cong K_2$ be the
restriction of $G$ on the edge $\overline{xy}\in E(G)$.

\medskip
Following the definitions from Section~\ref{sec:natural} the map
\[
\mathcal{F}^H : E(G) \longrightarrow Top,
\]
where $\mathcal{F}^H(\overline{xy}):=
\mathcal{F}^H_{\overline{xy}} = Hom(G_{\overline{xy}}, H)$, can be
thought of as a ``bundle'' over the graph $G$, with
$\mathcal{F}^H_{\overline{xy}}=Hom(G_{\overline{xy}},H)$ in the
role of the ``fibre'' over the edge $\overline{xy}$. More
generally, given a class $\mathcal{C}$ of subgraphs of $G$, say
the subtrees, the chains, the $k$-cliques etc., one can define an
associated ``bundle'' $\mathcal{F}^H_{\mathcal{C}}:
\mathcal{C}\rightarrow Top$ by a similar formula
$\mathcal{F}^H_{\mathcal{C}}(\Gamma):= Hom(\Gamma, H)$, where
$\Gamma\in \mathcal{C}$.

The parallel transport $\mathcal{P}^H$, for a given graph
($1$-dimensional, simplicial complex) $H$, is a specialization of
the parallel transport $\mathcal{P}^L$ introduced in
Section~\ref{sec:natural}. For example if
$\overrightarrow{e_1e_2}$ is the perspectivity between adjacent
edges $e_1=\overline{x_0x_1}$ and $e_2=\overline{x_1x_2}$ in $G$,
and if $\eta : \{x_1,x_2\}\rightarrow
2^{V(H)}\setminus\{\emptyset\}$ is a cell in
$\mathcal{F}^H_{\overline{x_1x_2}} = Hom(G_{\overline{x_1x_2}},
H)$, then $\eta':= \mathcal{P}^H(\overrightarrow{e_1e_2})(\eta) :
\{x_0,x_1\}\rightarrow 2^{V(H)}\setminus\{\emptyset\}$ is defined
by
\[
\eta'(x_0) := \eta(x_2)  \mbox{ {\rm and} } \eta'(x_1) :=
\eta(x_1).
\]

\medskip\noindent
{\bf Fundamental observation:} The construction of the connections
$\mathcal{P}^L$, respectively $\mathcal{P}^H$, are quite natural
and elementary but it is Proposition~\ref{prop:transport},
respectively its more general relative
Proposition~\ref{prop:prenos}, that serve as an actual
justification for the introduction of these objects.
Proposition~\ref{prop:transport} allows us to analyze the parallel
transport of homotopy types of maps from the complex $Hom(G,H)$ to
complexes $Hom(G_e,H)$, where $e\in E(G)$, providing a key for a
resolution of the Lov\'{a}sz conjecture in the case when $k$ is an
odd integer.

\medskip\noindent
Implicit in the proof of Proposition~\ref{prop:transport} is the
theory of {\em folds} of graphs and the analysis of natural
morphisms between graph complexes $Hom(T,H)$, where $T$ is a tree,
as developed in \cite{BabsonKozlov1} \cite{Kozl2004}
\cite{Kozlov99}. This theory is one of essential ingredients in
the Babson and Kozlov spectral sequence approach to the solution
of Lov\'{a}sz conjecture. Some of these results are summarized in
Proposition~\ref{prop:help}, in the form suitable for application
to Proposition~\ref{prop:transport}.

As usual $L_m$ is the graph-chain of vertex-length $m$, while
$L_{x_1\ldots x_m}$ is the graph isomorphic to $L_m$ defined on a
linearly ordered set of vertices $x_1,\ldots,x_m$. In this context
the ``flip'' is a generic name for the automorphism $\sigma :
L_{x_1\ldots x_m}\rightarrow L_{x_1\ldots x_m}$ of the graph-chain
such that $\sigma(x_j)=x_{m-j+1}$ for each $j$.

\begin{prop}\label{prop:help}
Suppose that $e_1=\overline{x_0x_1}$ and $e_2=\overline{x_1x_2}$
are two distinct, adjacent edges in the graph $G$. Let $\sigma :
L_{x_0x_1x_2}\rightarrow L_{x_0x_1x_2}$ be the flip automorphism
of $L_{x_0x_1x_2}$ and $\widehat{\sigma}$ the associated
auto-homeomorphism of $Hom(L_{x_0x_1x_2},H)$. Suppose that
$\gamma_{ij}: L_{x_ix_j}\rightarrow L_{x_0x_1x_2}$ is an obvious
embedding and $\widehat{\gamma}_{ij}$ the associated maps of graph
complexes. Then,
 \begin{enumerate}
 \item[{\rm (a)}] the induced map $\widehat{\sigma}:
Hom(L_{x_0x_1x_2},H)\rightarrow Hom(L_{x_0x_1x_2},H)$ is homotopic
to the identity map $I$, and
 \item[{\rm (b)}]  the diagram
\[
\begin{CD}
Hom(L_{x_0x_1x_2},H) @>=>> Hom(L_{x_0x_1x_2},H) \\
 @V\widehat{\gamma}_{01}VV @VV\widehat{\gamma}_{12}V \\
 Hom(L_{x_0x_1},H) @<<\mathcal{P}^H(\overrightarrow{e_1e_2})<
 Hom(L_{x_1x_2},H)
 \end{CD}
\]
 is commutative up to homotopy.
 \end{enumerate}
\end{prop}

\medskip\noindent
{\bf Proof:} Both statements are corollaries of Babson and Kozlov
analysis of complexes $Hom(T,H)$, where $T$ is a tree, and
morphisms $\widehat{e}: Hom(T,H)\rightarrow Hom(T',H)$, where $T'$
is a subtree of $T$ and $e: T'\rightarrow T$ the associated
embedding.

Our starting point is an observation that both $L_{x_0x_1}$ and
$L_{x_1x_2}$ are retracts of the graph $L_{x_0x_1x_2}$ in the
category of graphs and graph homomorphisms. The retraction
homomorphisms $\phi_{ij}: L_{x_0x_1x_2}\rightarrow L_{x_ix_j}$,
where $\phi_{01}(x_0)=x_0, \phi_{01}(x_1)=x_1, \phi_{01}(x_2)=x_0$
and $\phi_{12}(x_0)=x_2, \phi_{12}(x_1)=x_1, \phi_{12}(x_2)=x_2$
are examples of {\em foldings} of graphs. By the general theory
\cite{BabsonKozlov1} \cite{Kozl2004}, the maps
$\widehat{\gamma}_{ij}: Hom(L_{x_0x_1x_2},H)\rightarrow
Hom(L_{x_ix_j},H)$ and $\widehat{\phi}_{ij}:
Hom(L_{x_ix_j},H)\rightarrow Hom(L_{x_0x_1x_2},H)$ are homotopy
equivalences. Actually $\widehat{\gamma}_{ij}$ is a deformation
retraction and $\widehat{\phi}_{ij}$ is the associated embedding
such that $\widehat{\gamma}_{ij}\circ \widehat{\phi}_{ij}=I$ is
the identity map.

The part (a) of the proposition is an immediate consequence of the
fact that $\phi_{01}\circ\sigma\circ\gamma_{01} :
L_{x_0x_1}\rightarrow L_{x_0x_1}$ is an identity map. It follows
that
$\widehat{\gamma_{01}}\circ\widehat{\sigma}\circ\widehat{\phi}_{01}=I$,
and in light of the fact that $\widehat{\gamma_{01}}$ and
$\widehat{\phi}_{01}$ are homotopy inverses to each other, we
conclude that $\widehat{\sigma}\simeq I$.

For the part (b) we begin by an observation that $\phi_{12}\circ
\sigma\circ \gamma_{01} = \overrightarrow{e_1e_2}$. Then,
$\mathcal{P}^H(\overrightarrow{e_1e_2}) =
\widehat{\gamma}_{01}\circ\widehat{\sigma}\circ\widehat{\phi}_{12}$,
and as a consequence of $\widehat{\sigma} \simeq I$ and the fact
that $\widehat{\phi}_{12}\circ \widehat{\gamma}_{12}\simeq I$, we
conclude that
\[
\mathcal{P}^H(\overrightarrow{e_1e_2})\circ \widehat{\gamma}_{12}
= \widehat{\gamma}_{01}\circ\widehat{\sigma}\circ
\widehat{\phi}_{12}\circ\widehat{\gamma}_{12}\simeq
\widehat{\gamma}_{01}.
\]
\hfill $\square$

\begin{prop}\label{prop:transport} Suppose that $x_0,x_1,x_2$ are
distinct vertices in $G$ such that $\overline{x_0x_1},
\overline{x_1x_2}\in E(G)$. Let $\alpha_{ij}:
G_{x_ix_{j}}\rightarrow G$ be the inclusion map of graphs and
$\widehat{\alpha}_{ij}$ the associated map of $Hom(\,\cdot\,, H)$
complexes. Then the following diagram commutes up to a homotopy,

\begin{equation}\label{eqn:diagram1}
\begin{CD}
Hom(G,H) @>=>> Hom(G,H) \\
 @V\widehat{\alpha}_{01}VV @VV\widehat{\alpha}_{12}V \\
 Hom(G_{x_0x_1},H) @<<\mathcal{P}^H(\overrightarrow{e_1e_2})<
 Hom(G_{x_1x_2},H)
 \end{CD}
\end{equation}
\end{prop}

\medskip\noindent
{\bf Proof:} The diagram (\ref{eqn:diagram1}) can be factored as
\begin{equation}\label{eqn:diagram2}
\begin{CD}
  Hom(G,H) @>=>> Hom(G,H) \\
  @V\widehat{\beta} VV @VV\widehat{\beta} V \\
  Hom(G_{x_0x_1x_2},H) @>=>> Hom(G_{x_0x_1x_2},H)\\
  @V\widehat{\gamma}_{01}VV @VV\widehat{\gamma}_{12}V\\
  Hom(G_{x_0x_1},H) @<<\mathcal{P}^H(\overrightarrow{e_1e_2})<
  Hom(G_{x_1x_2},H)
 \end{CD}
\end{equation}
where $\beta$ and $\gamma_{ij}$ are obvious inclusions of
indicated graphs such that $\alpha_{ij}=\beta\circ\gamma_{ij}$.
Then the result is a direct consequence of
Proposition~\ref{prop:help}, part (b).
 \hfill $\square$

\section{Lov\'{asz-Babson-Kozlov result for odd $k$}}

The proof \cite{BabsonKozlov2} of Lov\' asz conjecture splits into
two main branches, corresponding to the parity of a parameter $n$,
where $n$ is an integer which enters the stage as the size of the
vertex set of the complete graph $K_n$.

The first branch relies on Theorem~2.3.\ (loc.\ cit.), more
precisely on part (b) of this result, while the second branch is
founded on Theorem~2.6. Both theorems are about the topology of
the graph complex $Hom(C_{2r+1}, K_n)$. Theorem~2.3.~(b) is a
statement about the height of the first Stiefel-Whitney class, or
equivalently the Conner-Floyd $\mathbb{Z}_2$-index
\cite{ConnerFloyd} of the $\mathbb{Z}_2$-space $Hom(C_{2r+1},
K_n)$. Theorem~2.6. claims that for $n$ even, $2\iota^\ast_{K_n}$
is a zero homomorphism where
\begin{equation}\label{eqn:thm-2.6.}
 \iota^\ast_{K_n} : \widetilde{H}^\ast(Hom(K_2,K_n);\mathbb{Z})
 \longrightarrow \widetilde{H}^\ast(Hom(C_{2r+1},K_n);\mathbb{Z})
\end{equation}
is the homomorphism  associated to the continuous map $\iota_{K_n}
: Hom(C_{2r+1},K_n)\rightarrow Hom(K_2,K_n)$, which in turn comes
from the inclusion $K_2\hookrightarrow C_{2r+1}$.

\medskip
The central idea of our paper is an observation that Theorem~2.6.\
can be incorporated into a more general scheme, involving the
``parallel transport'' of graph complexes over graphs.

\begin{thm}\label{thm:glavna}
Suppose that $\alpha : K_2 \rightarrow C_{2r+1}$ is an inclusion
map, $\beta : K_2 \rightarrow K_2$ a nontrivial automorphism of
$K_2$, and
 \[ \widehat{\alpha} : Hom(C_{2r+1},H)\rightarrow
Hom(K_2,H), \, \widehat{\beta} : Hom(K_2,H)\rightarrow Hom(K_2,H)
 \]
the associated maps of graph complexes. Then the following diagram
is commutative up to a homotopy
\begin{equation}\label{eqn:diagram3}
\begin{CD}
  Hom(C_{2r+1},H) @>=>> Hom(C_{2r+1},H) \\
  @V\widehat{\alpha} VV @VV\widehat{\alpha} V \\
  Hom(K_2,H) @<\widehat{\beta} << Hom(K_2,H)
   \end{CD}
\end{equation}
\end{thm}

\medskip\noindent
{\bf Proof:} Assume that the consecutive vertices of $G=C_{2r+1}$
are $x_0, x_1, \ldots , x_{2r}$ and let $e_i =
\overline{x_{i-1}x_i}$ be the associated sequence of edges where
by convention $e_{2r+1}=\overline{x_{2r}x_0}$. Identify the graph
$K_2$ to the subgraph $G_{x_0x_1}$ of $G=C_{2r+1}$.

By iterating Proposition~\ref{prop:transport} we observe that the
diagram
\begin{equation}\label{eqn:diagram4}
\begin{CD}
Hom(C_{2r+1},H) @>=>> Hom(C_{2r+1},H) \\
 @V\widehat{\alpha} VV @VV\widehat{\alpha} V \\
 Hom(G_{x_0x_1},H) @<<\mathcal{P}^H(\mathfrak{p})<
 Hom(G_{x_0x_1},H)
 \end{CD}
\end{equation}
is commutative up to a homotopy, where $\mathfrak{p}=
\overrightarrow{e_1e_2}\ast\ldots\ast\overrightarrow{e_{2r+1}e_1}$.
The proof is completed by the observation that $\mathfrak{p}=
\beta$ in the groupoid $\mathcal{G}^P(G)$. \hfill $\square$

\medskip
Theorem~2.6. from \cite{BabsonKozlov2}, the key for the proof of
Lov\'{a}sz conjecture for odd $k$, is an immediate consequence of
Theorem~\ref{thm:glavna}.

\begin{cor}{\rm (\cite{BabsonKozlov2}, T.2.6.)}
If $n$ is even then $2\cdot\iota^\ast_{K_n}$ is a $0$-map where
$\iota^\ast_{K_n}$ is the map described in line\ {\rm
(\ref{eqn:thm-2.6.})}.
\end{cor}

\medskip\noindent
{\bf Proof:} It is sufficient to observe that for $H=K_n$, the
complex $Hom(K_2,K_n)\cong S^{n-2}$ is an even dimensional sphere
such that the automorphism $\widehat{\beta}$ from the diagram
(\ref{eqn:diagram3}) is essentially an antipodal map. It follows
that $\widehat{\beta}$ changes the orientation of $Hom(K_2,K_n)$
and as a consequence $\iota^\ast_{K_n}= -\iota^\ast_{K_n}$. \hfill
$\square$

\section{Generalizations and ramifications}
\label{sec:ramifications}

In this section we extend the results from
Section~\ref{sec:functor} to the case of simplicial complexes.
This generalization is based on the following basic principles.

Graphs are viewed as $1$-dimensional simplicial complexes. Graph
homomorphisms are special cases of {\em non-degenerated}
simplicial maps of simplicial complexes, \cite{IzmJos2002}
\cite{Josw2001}. The definition of $Hom(G,H)$ is extended to the
case of $Hom$-complexes $Hom(K,L)$ of simplicial complexes $K$ and
$L$. The groupoids needed for the definition of the parallel
transport of $Hom$-complexes are already introduced by Joswig in
\cite{Josw2001}, see Section~\ref{sec:natural} for a summary.
Theory of folds for graph complexes \cite{BabsonKozlov1}
\cite{Kozl2004} is extended in Section~\ref{sec:tree-like} to the
case of $Hom$-complexes in sufficient generality to allow
``parallel transport'' of homotopy types of maps between graph
complexes. This development eventually leads to
Theorem~\ref{thm:main} which extends Theorem~\ref{thm:glavna} to
the case of $Hom$-complexes $Hom(K,L)$ and represents the
currently final stage in the evolution of Theorem~2.6. from
\cite{BabsonKozlov2}.

\begin{table}[hbt]\label{tab:tabela}
\begin{center}
\begin{tabular}{||l|l||}\hline
\multicolumn{2}{||c||}{\bf Dictionary}\\ \hline\hline
 graphs &  simplicial complexes \\ \hline
 trees & tree-like complexes \\ \hline
 foldings of graphs & vertex collapsing of complexes \\ \hline
 graph homomorphisms & non-degenerated simplicial maps \\ \hline
 $Hom(G,H)$ & $Hom(K,L)$ \\ \hline
 chromatic number $\chi(G)$ & chromatic number $\chi(K)$ \\ \hline
\end{tabular}
\end{center}
\caption{Graphs vs.\ simplicial complexes.}
\end{table}

\subsection{From $Hom(G,H)$ to $Hom(K,L)$}
\label{sec:Hom}

Suppose that $K\subset 2^{V(K)}$ and $L\subset 2^{V(L)}$ are two
(finite) simplicial complexes, on the sets of vertices $V(K)$ and
$V(L)$ respectively.

\begin{defn}
A simplicial map $f: K\rightarrow L$ is {\em non-degenerated} if
it is injective on simplices. The set of all non-degenerated
simplicial maps from $K$ to $L$ is denoted by $Hom_0(K,L)$.
\end{defn}

\begin{defn}\label{def:hom(k,l)}
$Hom(K,L)$ is a cell complex with the cells indexed by the
functions $\eta : V(K)\rightarrow 2^{V(L)}\setminus\{\emptyset\}$
such that
 \begin{enumerate}
 \item[{\rm (1)}] for each two vertices $u\neq v$, if $\{u,v\}\in K$ then
$\eta(u)\cap\eta(v)=\emptyset$,
 \item[{\rm (2)}] for each simplex $\sigma\in K$, the join
 ${{\ast}}_{v\in V(\sigma)}~\eta(v)\subset \Delta^{V(L)}$
 of all {\em sets ($0$-dimensional complexes)} $\eta(v)$, where $v$ is a
 vertex of $\sigma$, is a subcomplex of $L$.
\end{enumerate}
More precisely, each function $\eta$ satisfying conditions {\rm
(1)} and {\rm (2)} defines a cell $c_\eta:= \prod_{v\in
V(K)}~\Delta^{\eta(v)}$ in $Hom(K,L)\subset\prod_{v\in
V(K)}\Delta^{V(L)}$ where by definition $\Delta^S$ is an
(abstract) simplex spanned by vertices in $S$.
\end{defn}

We have already used in Section~\ref{sec:natural} the fact that if
$K=\Delta^{[m]}$ is a $(m-1)$-dimensional simplex spanned by $[m]$
as the set of vertices, then $Hom(\Delta^{[m]},L)\cong
L_\Delta^{m}$ is the {\em deleted product} of $L$ \cite{Matousek}.
The following example shows that $Hom(G,H)$ is a special case of
$Hom(K,L)$.

\begin{exmp}{\rm The definition of $Hom(K,L)$ is a natural
extension of $Hom(G,H)$ and reduces to it if $K$ and $L$ are
$1$-dimensional complexes. Moreover,
\[
Hom(G,H)\cong Hom(Clique(G),Clique(H))
\]
where $Clique(\Gamma)$ is the simplicial complex of all cliques in
a graph $\Gamma$. }
\end{exmp}

\begin{rem}{\rm
The set $Hom_0(K,L)$ is easily identified as the $0$-dimensional
skeleton of the cell-complex $Hom(K,L)$. Moreover, the reader
familiar with \cite{Kozlov-review} can easily check that
$Hom(K,L)$ is determined by the family $M=Hom_0(K,L)$ in the sense
of Definition~2.2.1. from that paper. }
\end{rem}

\subsection{Functoriality of $Hom(K,L)$}
\label{sec:naturality}

The construction of $Hom(K,L)$ is functorial in the sense that if
$f : K\rightarrow K'$ is a non-degenerated simplicial map of
complexes $K$ and $K'$, then there is an associated continuous map
$\widehat{f} : Hom(K',L)\rightarrow Hom(K,L)$ of $Hom$-complexes.
Indeed, if $\eta : V(K')\rightarrow
2^{V(L)}\setminus\{\emptyset\}$ is a multi-valued function
indexing a cell in $Hom(K',L)$, then it is not difficult to check
that $\eta\circ f : V(K)\rightarrow
2^{V(L)}\setminus\{\emptyset\}$ is a cell in $Hom(K,L)$.

Perhaps even more important is the functoriality of $Hom(K,L)$
with respect to the second variable since this implies the
functoriality of the bundle $\mathcal{F}^L_k$.

\begin{prop}\label{prop:functoriality}
Suppose that $g : L \rightarrow L'$ is a non-degenerated,
simplicial map of simplicial complexes $L$ and $L'$. Then there
exists an associated map
\[
\widehat{g} : Hom(K,L)\rightarrow Hom(K,L').
\]
\end{prop}

\medskip\noindent
{\bf Proof:} Assume that $\eta : V(K)\rightarrow
2^{V(L)\setminus\{\emptyset\}}$ is a cell in $Hom(K,L)$. Then
$g\circ\eta : V(K)\rightarrow 2^{V(L')\setminus\{\emptyset\}}$ is
a cell in $Hom(K,L')$. Suppose $u$ and $v$ are distinct vertices
in $V(K)$. By assumption $\eta(u)\cap \eta(v)=\emptyset$. We
deduce from here that $g(\eta(u))\cap g(\eta(v))\neq \emptyset$,
otherwise $g$ would be a degenerated simplicial map.

The second condition from Definition~\ref{def:hom(k,l)} is checked
by a similar argument. \hfill $\square$

\subsection{Chromatic number $\chi(K)$ and its relatives}

The chromatic number $\chi(K)$ of a simplicial complex $K$ is
\[
\inf\{m\in \mathbb{N}\mid Hom_0(K,\Delta^{[m]})\neq\emptyset\}.
\]
In other words $\chi(K)$ is the minimum number $m$ such that there
exists a non-degenerated simplicial map $f : K\rightarrow
\Delta^{[m]}$. It is not difficult to check that
$\chi(K)=\chi(G_K)$ where $G_K = (K^{(0)}, K^{(1)})$ is the
vertex-edge graph of the complex $K$. In particular $\chi(K)$
reduces to the usual chromatic number if $K$ is a graph, that is
if $K$ is a $1$-dimensional simplicial complex.

\medskip
Aside from the usual chromatic  number $\chi(G)$, there are many
related colorful graph invariants \cite{GodRoy}
\cite{Kozlov-review}. Among the best known are the fractional
chromatic number $\chi_f(G)$ and the circular chromatic number
$\chi_c(G)$ of $G$. These and other related invariants are
conveniently defined in terms of graph homomorphisms into graphs
chosen from a suitable family $\mathcal{F} =\{G_i\}_{i\in I}$ of
test graphs. Motivated by this, we offer an extension of the
chromatic number $\chi(K)$ in hope that some genuine invariants of
simplicial complexes objects arise this way.

\begin{defn}
Suppose that $\mathcal{F}=\{T_i\mid i\in I\}$ is a family of
``test'' simplicial complexes and let $\phi : I\rightarrow
\mathbb{R}$ is a real-valued function. A $T_i$-coloring of $K$ is
just a non-degenerated simplicial map $f:K\rightarrow T_i$ and
$\chi_{(\mathcal{F},\phi)}(K)$, the
$(\mathcal{F},\phi)$-chromatic number
 of $K$, is defined as the infimum
of all weights $\phi(i)$ over all $T_i$-colorings,

\[
\chi_{(\mathcal{F},\phi)}(K):= \inf\{\phi(i) \mid
Hom_0(K,T_i)\neq\emptyset\}.
\]

\end{defn}

\subsection{Tree-like simplicial complexes}
\label{sec:tree-like}

The tree-like or vertex collapsible complexes are intended to play
in the theory of $Hom(K,L)$-complexes the role similar to the role
of trees in the theory of graph complexes $Hom(G,H)$.

A pure, $d$-dimensional simplicial complex $K$ is {\em shellable}
\cite{Bjo-et-al} \cite{Zieg-book}, if there is a linear order
$F_1, F_2, \ldots , F_m$ on the set of its facets, such that for
each $j\geq 2$, the complex $F_j\cap(\bigcup_{i<j}~F_i)$ is a pure
$(d-1)$-dimensional subcomplex of the simplex $F_j$. The {\em
restriction} $R_j$ of the facet $F_j$ is the minimal new face
added to the complex $\bigcup_{k<j}~F_k$ by the addition of the
facet $F_j$. Let $r_j := {\rm dim}(R_j)\in\{0,1,\ldots,d\}$ be the
type of the facet $F_j$. If $r_j\neq d$ for each $j$ then the
complex $K$ is collapsible. The collapsing process is just the
shelling order read in the opposite direction. From this point of
view, $R_j$ can be described as a free face in the complex
$\bigcup_{i\leq j}~F_i$, and the process of removing all faces $F$
such that $R_j\subset F\subset F_j$ is called an elementary
$r_j$-collapse.

\begin{defn}
A pure $d$-dimensional simplicial complex $K$ is called {\em
tree-like} or vertex collapsible if it is collapsible to a
$d$-simplex with the use of elementary $0$-collapses alone. In
other words $K$ is shellable and for each $j\geq 2$, the
intersection $F_j\cap(\bigcup_{i<j}~F_i)$ is a proper face of
$F_j$.
\end{defn}

In order to establish analogs of Propositions~\ref{prop:help} and
\ref{prop:transport} for complexes $Hom(K,L)$, we prove a result
which shows that elementary vertex collapsing  provides a good
substitute and a partial generalization for the concept of
``foldings'' of graphs used in \cite{BabsonKozlov1}
\cite{Kozl2004} in the theory of graph complexes $Hom(G,H)$.

\begin{prop}
\label{prop:vertex-collapse} Suppose that the simplicial complex
$K'$ is obtained from $K$ by an elementary vertex collapse. In
other words we assume that $K = \sigma\cup K',\, \sigma\cap
K'=\sigma'$, where $\sigma$ is a simplex in $K$ and $\sigma'$ a
facet of $\sigma$. Assume that $\sigma'$ is not maximal in $K'$,
i.e.\ that for some simplex $\sigma''\in K'$ and a vertex
$u\in\sigma'',\, \sigma' = \sigma''\setminus\{u\}$. Then for any
simplicial complex $L$, the inclusion map $\gamma : K'\rightarrow
K$ induces a homotopy equivalence
\[
\widehat{\gamma} : Hom(K,L)\rightarrow Hom(K',L).
\]
\end{prop}

\medskip\noindent
{\bf Proof:} Let $\{v\}=\sigma\setminus \sigma'$. Aside from the
inclusion map $\gamma : K'\rightarrow K$, there is a retraction
(folding) map $\rho : K\rightarrow K'$, where $\rho(v)=u$ and
$\rho\vert_{K'}=I_{K'}$. Since $\rho\circ\gamma = I_{K'}$, we
observe that $\widehat{\gamma}\circ\widehat{\rho}=Id_{K'}$ is the
identity map on $Hom(K',L)$, i.e.\ the complex $Hom(K',L)$ is a
retract of the complex $Hom(K,L)$. It remains to be shown that
$\widehat{\rho}\circ\widehat{\gamma}\simeq Id_{K}$ is homotopic to
the identity map on $Hom(K,L)$.

Note that if $\eta\in Hom(K,L)$ then $\eta':=
\widehat{\rho}\circ\widehat{\gamma}(\eta)$ is the function defined
by
 $$
 \eta'(w) = \left\{\begin{array}{rl} \eta(w), & \mbox{ {\rm
if} \, }
w\neq v \\
\eta(u), & \mbox{ {\rm if} \, } w = v.
\end{array}\right.
 $$
Let $\omega : Hom(K,L)\rightarrow Hom(K,L)$ be the map defined by
$$
 \omega(\eta)(w) = \left\{\begin{array}{ll} \eta(w), & \mbox{ {\rm
if} \, }
w\neq v \\
\eta(u)\cup\eta(v), & \mbox{ {\rm if} \, } w = v.
\end{array}\right.
 $$
Note that $\omega$ is well defined since if a vertex $x$ is
adjacent to $v$ it is also adjacent to $u$, hence the condition
$\omega(\eta)(v)\cap\eta(x)=\emptyset$ is a consequence of
$\eta(u)\cap\eta(x)=\emptyset = \eta(v)\cap\eta(x)$.

Since for each $\eta\in Hom(K,L)$ and each vertex $x\in K$,
\[
\eta(x)\subset \omega(\eta)(x) \supset
\widehat{\rho}\circ\widehat{\gamma}(\eta)(x),
\]
by the Order Homotopy Theorem \cite{Bjorner} \cite{Quillen}
\cite{Segal} all three maps $Id_{K}, \omega$ and
$\widehat{\rho}\circ\widehat{\gamma}$ are homotopic. This
completes the proof of the proposition. \hfill $\square$

\begin{cor}
If $T$ is a $d$-dimensional, tree-like simplicial complex than
$Hom(T,L)$ has the same homotopy type as the deleted join
$Hom(\Delta^d,L)=L^{d+1}_\Delta$.
\end{cor}

\subsection{Parallel transport of homotopy types of maps}

As in the case of graph complexes, the real justification for the
introduction of the parallel transport of $Hom$-complexes comes
from the fact that it preserves the homotopy type of the maps
$Hom(K,L)\rightarrow Hom(\sigma, L)$. As in
Section~\ref{sec:functor}, as a preliminary step we prove an
analogue of Proposition~\ref{prop:help}.

\begin{prop}\label{prop:helphelp}
Suppose that $\sigma_1$ and $\sigma_2$ are two distinct, adjacent
$k$-dimensional simplices in a finite simplicial complex $K$ which
share a common $(k-1)$-dimensional simplex $\tau$. Let $\Sigma =
\sigma_1\cup\sigma_2$. Let $\alpha : \Sigma\rightarrow\Sigma$ be
the automorphism of $\Sigma$ which interchanges simplices
$\sigma_1$ and $\sigma_2$ keeping the common face $\tau$
point-wise fixed.

Suppose that $\gamma_i: \sigma_i \rightarrow \Sigma$ is an obvious
embedding and $\widehat{\gamma}_i$ the associated maps of
$Hom$-complexes. Then,
 \begin{enumerate}
 \item[{\rm (a)}] the induced map $\widehat{\alpha}:
Hom(\Sigma,L)\rightarrow Hom(\Sigma,L)$ is homotopic to the
identity map $I_\Sigma$, and
 \item[{\rm (b)}]  the diagram
\[
\begin{CD}
Hom(\Sigma,L) @>=>> Hom(\Sigma,L) \\
 @V\widehat{\gamma}_{1}VV @VV\widehat{\gamma}_{2}V \\
 Hom(\sigma_1,L) @<<\mathcal{P}^H(\overrightarrow{\sigma_1\sigma_2})<
 Hom(\sigma_2,L)
 \end{CD}
\]
is commutative up to homotopy.
 \end{enumerate}
\end{prop}

\medskip\noindent
{\bf Proof:} By Proposition~\ref{prop:vertex-collapse}, both maps
$\widehat{\gamma}_i : Hom(\Sigma,L)\rightarrow Hom(\sigma_i,L)$
for $i=1,2$ are homotopy equivalences. Let $\rho_1 :
\Sigma\rightarrow \sigma_1$ and $\rho_2 : \Sigma\rightarrow
\sigma_2$ be the folding maps. Then $\rho_i\circ\gamma_i =
I_{\sigma_i},\, $ $\widehat{\gamma_i}\circ\widehat{\rho_i}= I$ and
we conclude that $\widehat{\rho}_i : Hom(\sigma_i,L)\rightarrow
Hom(\Sigma,L)$ is also a homotopy equivalence.

Part (a) of the proposition follows from the fact that
$\rho_1\circ\alpha\circ\gamma_1 = I_{\sigma_1}$ is an identity
map. Indeed, an immediate consequence is that
$\widehat{\gamma_1}\circ \widehat{\alpha}\circ \widehat{\rho_1} =
I : Hom(\sigma_1,L)\rightarrow Hom(\sigma_1,L)$ is also an
identity map and, in light of the fact that $\widehat{\gamma_1}$
and $\widehat{\rho_1}$ are homotopy inverses to each other, we
deduce that $\widehat{\alpha}\simeq I$.

For the part (b) we begin by an observation that $\rho_2\circ
\alpha\circ \gamma_1 = \overrightarrow{\sigma_1\sigma_2}$. Then,
$\mathcal{P}^H(\overrightarrow{\sigma_1\sigma_2}) =
\widehat{\gamma}_1\circ\widehat{\alpha}\circ\widehat{\rho}_2$, and
as a consequence of $\widehat{\alpha} \simeq I$ and the fact that
$\widehat{\rho}_2\circ \widehat{\gamma}_2\simeq I$, we conclude
that
\[
\mathcal{P}^H(\overrightarrow{\sigma_1\sigma_2})\circ
\widehat{\gamma}_2 = \widehat{\gamma}_1\circ\widehat{\alpha}\circ
\widehat{\rho}_2\circ\widehat{\gamma}_2\simeq \widehat{\gamma}_1.
\]
\hfill $\square$

\begin{prop}\label{prop:prenos}
Suppose that $K$ and $L$ are finite simplicial complexes and
$\sigma_1,\sigma_2$ a pair of adjacent (distinct), $k$-dimensional
simplices in $K$. Let $\alpha_i : \sigma_i\rightarrow K$ be the
embedding of $\sigma_i$ in $K$ and $\widehat{\alpha}_i:
Hom(K,L)\rightarrow Hom(\sigma_i,L)$ the associated map of
$Hom$-complexes. Then the following diagram commutes up to a
homotopy.

\begin{equation}\label{eqn:shema1}
\begin{CD}
Hom(K,L) @>=>> Hom(K,L) \\
 @V\widehat{\alpha}_{1}VV @VV\widehat{\alpha}_{2}V \\
 Hom(\sigma_1,L) @<<\mathcal{P}^L(\overrightarrow{e_1e_2})<
 Hom(\sigma_2,H)
 \end{CD}
\end{equation}
\end{prop}

\medskip\noindent
{\bf Proof:} Let $\Sigma := \sigma_1\cup\sigma_2, \, \tau :=
\sigma_1\cap\sigma_2$. Then $\alpha_i = \beta\circ\gamma_i$ where
$\beta:\Sigma\rightarrow K$ and $\gamma_i: \sigma_i \rightarrow
\Sigma$ are natural embeddings of complexes.

The diagram (\ref{eqn:shema1}) can be factored as
\begin{equation}\label{eqn:shema2}
\begin{CD}
  Hom(K,L) @>=>> Hom(K,L) \\
  @V\widehat{\beta} VV @VV\widehat{\beta} V \\
  Hom(\Sigma,L) @>=>> Hom(\Sigma,L)\\
  @V\widehat{\gamma}_1VV @VV\widehat{\gamma}_2 V\\
  Hom(\sigma_1,L) @<<\mathcal{P}^L(\overrightarrow{\sigma_1\sigma_2})<
  Hom(\sigma_2,L)
 \end{CD}
\end{equation}
Then the result is a direct consequence of
Proposition~\ref{prop:helphelp}, part (b).
 \hfill $\square$

\begin{cor}\label{cor:posledica}
Suppose that $K$ and $L$ are finite simplicial complexes, $\sigma$
a $k$-dimensional simplex in $K$ and $\alpha : \sigma\rightarrow
K$ the associated embedding. Let $\tau\in \Pi(K,\sigma)$. Then the
following diagram commutes up to a homotopy.

\begin{equation}\label{eqn:shema3}
\begin{CD}
Hom(K,L) @>=>> Hom(K,L) \\
 @V\widehat{\alpha} VV @VV\widehat{\alpha} V \\
 Hom(\sigma,L) @<<\widehat{\tau} <
 Hom(\sigma,L)
 \end{CD}
\end{equation}
\end{cor}

\section{Main results}
\label{sec:main}

In this section we prove the promised extension of the
Lov\'{a}sz-Babson-Kozlov theorem. The graphs are replaced by pure
$d$-dimensional simplicial complexes, while the role of the odd
cycle $C_{2r+1}$ is played by a complex $\Gamma$ which has some
special symmetry properties in the sense of the following
definition.

As usual, an involution $\omega : X\rightarrow X$ is the same as a
$\mathbb{Z}_2$-action on $X$. An involution on a simplicial
complex $\Gamma$ induces an involution on the complex
$Hom(\Gamma,L)$ for each simplicial complex $L$. For all other
standard facts and definitions related to
$\mathbb{Z}_2$-complexes, the reader is referred to
\cite{Matousek}.

\begin{defn}
\label{def:C-complex} A pure $d$-dimensional simplicial complex
$\Gamma$ is a $\Phi_d$-complex if it is a $\mathbb{Z}_2$-complex
with an invariant $d$-simplex $\sigma = \{v_0,v_1,\ldots, v_d\}$
such that the restriction $\tau :=\omega\vert_\sigma$ of the
involution $\omega : \Gamma\rightarrow\Gamma$ on $\sigma$ is a
non-trivial element of the group $\Pi(\Gamma,\sigma)$.
\end{defn}

\begin{rem}{\rm
By definition, if $\Gamma$ is a $\Phi_d$-complex then the
inclusion map $\alpha: \sigma\rightarrow\Gamma$ is
$\mathbb{Z}_2$-equivariant, so the associated map
$\widehat{\alpha}:Hom(\Gamma,K)\rightarrow Hom(\sigma,K)$ is also
$\mathbb{Z}_2$-equivariant for each complex $K$. }
\end{rem}

\begin{figure}[hbt]
\centering
\includegraphics[scale=0.40]{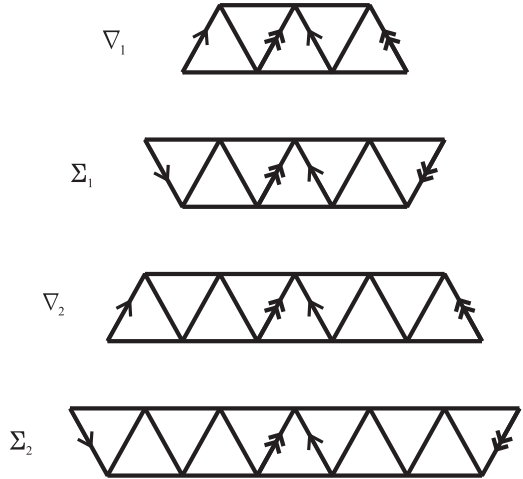}
\caption{Examples of $\Phi_2$-complexes.} \label{fig:piramida}
\end{figure}

\begin{exmp}{\rm
The graph $C_{2r+1}$ is obviously an example of a
$\Phi_1$-complex. Figure~\ref{fig:piramida} displays four examples
of $\Phi_2$-complexes, initial elements of two infinite series
$\nabla_\mu$ and $\Sigma_\nu,\, \mu,\nu\in \mathbb{N}$. The
complexes $\nabla_1$ and $\nabla_2$ etc.\ are obtained from two
triangulated annuli, glued together along a common triangle
$\sigma$. Similarly, the complexes $\Sigma_1, \Sigma_2, \ldots$,
are obtained by gluing together two triangulated M\"{o}bius
strips. The associated group of projectivities are
$\Pi(\nabla_\mu,\sigma)=S_3$ and
$\Pi(\Sigma_\nu,\sigma)=\mathbb{Z}_2$. }
\end{exmp}

\begin{thm}\label{thm:main}
Suppose that $\Gamma$ is a $\Phi_d$-complex in the sense of
Definition~\ref{def:C-complex}, with an associated invariant
simplex $\sigma = \{v_0,v_1,\ldots, v_d\}$. Suppose that $K$ is a
pure $d$-dimensional simplicial complex. Than for $m$ even,
\begin{equation}\label{eq:complex}
Coind_{\mathbb{Z}_2}(Hom(\Gamma,K))\geq m \Rightarrow \chi(K)\geq
m+d+2.
\end{equation}
\end{thm}

\medskip\noindent
{\bf Proof:} By definition
$Coind_{\mathbb{Z}_2}(Hom(\Gamma,K))\geq m$ means that there
exists a $\mathbb{Z}_2$-equiva\-ri\-ant map $\mu: S^m\rightarrow
Hom(\Gamma, K)$. Assume that $\chi(K)\leq m+d+1$ which means that
there exists a non-degenerated simplicial map $\phi :
K\rightarrow\Delta^{[m+d+1]}$. By functoriality of the
construction of $Hom$-complexes, Section~\ref{sec:naturality},
there is an induced $\mathbb{Z}_2$-equivariant map
$\widehat{\phi}: Hom(\Gamma, K)\rightarrow Hom(\Gamma,
\Delta^{[m+d+1]})$ and similarly a map $\widehat{\alpha}:
Hom(\Gamma, \Delta^{[m+d+1]})\rightarrow Hom(\sigma,
\Delta^{[m+d+1]})$. By \cite{Kozl2004} Theorem~3.3.3., the complex
 \[
Hom(\sigma, \Delta^{[m+d+1]})\cong Hom(K_{d+1},K_{m+d+1})
 \]
is a wedge of $m$-dimensional spheres. Since $Hom(\sigma,
\Delta^{[m+d+1]})$ is a free $\mathbb{Z}_2$-complex, we deduce
that there exists a $\mathbb{Z}_2$-equivariant  map $Hom(\sigma,
\Delta^{[m+d+1]})\rightarrow S^m$. All these maps can be arranged
in the following sequence of $\mathbb{Z}_2$-equivariant maps
 \[
S^m   \stackrel{\mu}\longrightarrow Hom(\Gamma, K)
\stackrel{\widehat{\phi}}\longrightarrow Hom(\Gamma,
\Delta^{[m+d+1]}) \stackrel{\widehat{\alpha}}\longrightarrow
Hom(\sigma, \Delta^{[m+d+1]}) \stackrel{\nu}\longrightarrow S^m.
 \]
By Corollary~\ref{cor:posledica}, there is a homotopy equivalence
$\widehat{\alpha}\simeq \tau\circ \widehat{\alpha}$. This is in
contradiction with Proposition~\ref{prop:not-homotopic}, which
completes the proof of the theorem. \hfill $\square$

\begin{prop}\label{prop:not-homotopic}
Suppose that $f: X\rightarrow Y$ is a $\mathbb{Z}_2$-equivariant
map of free $\mathbb{Z}_2$-complexes $X$ and $Y$ where
$\mathbb{Z}_2=\{1,\omega\}$. Assume that
$Coind_{\mathbb{Z}_2}(X)\geq m\geq Ind_{\mathbb{Z}_2}(Y)$, where
$m$ is an even integer. In other words our assumption is that
there exist $\mathbb{Z}_2$-equivariant maps $\mu$ and $\nu$ such
that
\[
S^m \stackrel{\mu}\longrightarrow X \stackrel{f}\longrightarrow Y
\stackrel{\nu}\longrightarrow S^m .
\]
Then the maps $f$ and $\omega\circ f$ are not homotopic.

\medskip\noindent
{\bf Proof:} If $f\simeq \omega\circ f : X\rightarrow Y$ then
$\nu\circ f\circ\mu \simeq \nu\circ\omega\circ f\circ\mu :
S^m\rightarrow S^m$ and by the equivariance of $\nu\,$,
$\omega\circ g\simeq g : S^m\rightarrow S^m$ where $g:= \nu\circ
f\circ \mu$. It follows that
\[
-{\rm deg}(g) = {\rm deg}(\omega) {\rm deg}(g) = {\rm
deg}(\omega\circ g) = {\rm deg}(g),
\]
i.e.\ ${\rm deg}(g)=0$, which is in contradiction with a well
known fact that a $\mathbb{Z}_2$-equivariant map $g :
S^m\rightarrow S^m$ of even dimensional spheres must have an odd
degree. \hfill $\square$
\end{prop}

\begin{cor}\label{cor:LBK-gen}
Suppose that $\Gamma$ is a $\Phi_d$-complex with an associated
invariant simplex $\sigma = \{v_0,v_1,\ldots, v_d\}$. Suppose that
$K$ is a pure $d$-dimensional simplicial complex. Than for $k$
odd,
\begin{equation}\label{eq:za-kraj}
Hom(\Gamma,K) \mbox{ {\rm is $k$-connected} } \Rightarrow
\chi(K)\geq k+d+3.
\end{equation}
\end{cor}

\medskip\noindent
{\bf Proof:} If $Hom(\Gamma,K)$ is $k$-connected then $Coind_{
\mathbb{Z}_2}(Hom(\Gamma,K))\geq k+1$, hence the implication
(\ref{eq:za-kraj}) is an immediate consequence of
Theorem~\ref{thm:main}. \hfill $\square$

\end{document}